\documentclass[2p,12pt,a4paper,twoside,fleqn,sort&compress]{amsart}
\usepackage{amssymb,amsmath,latexsym}
\usepackage[varg]{pxfonts}
\usepackage[english]{babel}
\newcommand{\be}{\begin{equation}}
\newcommand{\ee}{\end{equation}}

\newtheorem{thm}{Theorem}[section]
\newtheorem{pro}{Proposition}[section]
\newtheorem{lem}{Lemma}[section]
\newtheorem{rem}{Remark}[section]
\newtheorem{cor}{Corollary}[section]
\newtheorem{defin}{Definition}[section]

\DeclareFontFamily{T1}{cmr}{\hyphenchar\font=-1}

\newcommand{\pv}{\par \vskip0.2cm}

\numberwithin{equation}{section}
\DeclareFontFamily{T1}{cmr}{\hyphenchar\font=-1}
\usepackage{color}
\usepackage{latexsym}
\usepackage{cite}
\usepackage[numbers,sort&compress,square]{natbib}
\begin{document}
\title[]{Existence results for nonexpansive multi-valued
operators and
Nonlinear integral inclusions
}
\maketitle

\author\centerline { Khaled Ben Amara$^{(1)},$  Aref Jeribi $^{(1)}$ and   Najib Kaddachi $^{(2)}$ }
\date{}
\maketitle

\vskip0.3cm

\begin{center}
 Department of Mathematics. Faculty of Sciences of Sfax.
University of Sfax.\\ Road Soukra Km $3.5 B.P. 1171, 3000,$ Sfax,
Tunisia.\pv
\end{center}

\vskip0.3cm

\centerline {e-mail: $^{(1)}$ Aref.Jeribi$$@$$fss.rnu.tn,   \
$^{(2)}$ najibkadachi$$@$$gmail.com   }
\vskip0.7cm
\vskip0.2cm


\maketitle
\renewcommand{\thefootnote}{}

\noindent {\bf\small Abstract.}
In this paper, we establish some new variants of fixed point theorems for a large class of  countably nonexpansive multi-valued mappings. Some fixed point theorems for the sum and the product of three multi-valued mappings defined on nonempty, closed convex set of Banach algebras are also presented.
These results improve and complement a number of earlier
works.
 As an application, we prove existence results for a broad class of
nonlinear functional integral inclusions as well as nonlinear differential inclusions.

\vskip0.7cm
\par \noindent {\bf\small Keywords:}
 {\small \sloppy{Countably condensing multi-valued,  Countably nonexpansive multi-valued, Functional-differential inclusions, Banach algebras, Measure of weak noncompactness,  Fixed point theorems.
 }}
\vskip0.2cm
\par \noindent{\bf AMS Classification:}  $47$H$08,$ $47$H$09,$ $34$K$09,$ $47$B$48,$ $47$H$10.$
\vskip0.7cm
\section{Introduction}
Many problems arising in mathematical physics, chemistry, biology, medicine,
etc., can be described, in a first formulation, using nonlinear differential inclusions as well
as nonlinear integral inclusions, see \textup{\cite{agarwal, BA2016,Existence,Kaddachi2016,Kaddachi2018, JKM2013,Kaddachi2019,JK2015, Mitchell}}; for this reason, the development
of new tools allowing to solve such problems is of great interest.

Fixed point theorems for nonexpansive mappings is one of the important generalizations of the well-known Banach fixed point theorem.
 That is why several authors have focused on the existence
of fixed points of nonexpansive mappings in Banach spaces and have obtained many valuable results. We can cite for examples Ben Amar et al. \cite{Ben amar 16}, Browder \cite{N6}, Gohde \cite{N11}, Ishikawa \cite{N12}, Kirk \cite{N13}.
These results have been extended  by several authors to the case of multi-valued
 mappings. For examples, see \textup{\cite{agarwal,BADO2016,Husain80,Husain88}}.
In \cite{BADO2016}, A. Ben Amar and D. O'Regan have established some fixed point theorems for  $\beta$-nonexpansive multi-valued mappings, i.e.,
  \begin{eqnarray}\label{exp}\beta(T(M )) \leq\beta(M ),\end{eqnarray}
  for all bounded subset  $M$ of $\Omega$ such that  $T(M)$ is bounded.
 Here $\beta$ is the measure of weak noncompactness of De Blasi, see \cite{De Blasi}.

 In this paper, we investigate the existence of fixed point theorems for countably $\beta$-nonexpansive multi-valued mappings,  i.e., the condition (\ref{exp}) holds only for countable
bounded sets.
 Therefore, we prove new fixed point theorems of Krasnoselskii's type for multi-valued mappings defined on Banach spaces under weak topology. Moreover, we establish  the existence of solutions for the following hybrid fixed point inclusion on Banach algebras under the weak topology setting:
\begin{eqnarray}\label{in0}x \in Ax \cdot Bx+Cx.\end{eqnarray}
 The obtained results significantly extend and generalize some works in the literature.
 Our results are applied to discuss the existence of solutions
  to an abstract class of nonlinear integral inclusions of the form:
\begin{eqnarray}\label{in12}
   \displaystyle x(t)\in \displaystyle  T_1(t,x(t))  \cdot\left(q(t)+\int_0^{t}k(t,s)\,H(s,x(s))\, ds\right)+T_2(t,x(t)),\ \ \, t\in J,
\end{eqnarray}
 as well as the following class of nonlinear functional differential inclusions
\begin{eqnarray}\label{1}
\left\{
  \begin{array}{ll}
     \left(\displaystyle\frac{ x(t)-T_2(t,x(t))}{T_1(t,x(t))}\right)'\in k(t)\,H(t,x(t)); \\\\
    x(0)=x_0\in X,\ \ \, t\in J,
  \end{array}
\right.
\end{eqnarray}
  where $k : J \longrightarrow \mathbb{R},$ $q : J \longrightarrow X,$ $T_i : J \times X \longrightarrow X$ and $ H : J \times X \longrightarrow \mathcal{P}(X)$ with $X$ a Banach algebra.
\\
Some specials cases of the nonlinear inclusion (\ref{1}) have been studied in \textup{\cite{Dhage05I, Dhage05II,Dhage06}}.

This work is organized as follows, Section $2$ is devoted to some definitions and mainly to the
basic tools which will be used in the sequel. The section $3$
 deals with some fixed point results for countably $\mathcal{D}$-set-Lipschtzian and countably $\beta$-nonexpansive multi-valued mappings in Banach spaces under a weak topology setting. Moreover, we prove some fixed point theorems of Krasnoselskii's type for multi-valued mappings defined on Banach spaces.
 In section $4,$ we prove the existence of solutions for the inclusion (\ref{in0}).
 In the last section,  we
investigate an existence theory of solutions for the nonlinear integral inclusion $(\ref{in12})$ and the differential integral inclusion (\ref{1}).
\section{\textbf{Basics results}}
Let $E$ be a Banach space endowed with the norm $\|\cdot\|$ and with the zero element $\theta.$ For any $r>0,$ $\mathcal{B}_E(x,r) $ denotes the closed ball of $E$ centered at $x$ with radius $r,$ in particular $\mathcal{B}_E(r) :=\mathcal{B}_E(\theta,r),$
and for any subset $S$ of $E,$ we write $co(S)$ and $\overline{co}(S)$ to denote the convex hull and the
closed convex hull of $S,$ respectively.
The sets $\mathcal{B}(E)$
and $\mathcal{W}(E)$ stand, respectively, for the family of all nonempty
bounded subsets of $E$ and all nonempty weakly compact subsets of $E.$ 
 %
Let
\begin{eqnarray*}\displaystyle \mathcal{P}(E) &=& \big\{S \subset E:\, S\neq\emptyset\big\},\\
\displaystyle \mathcal{P}_{cv}(E)& = &\big\{S \in \mathcal{P}(E):\, S \hbox{ is convex}\big\},\\
\displaystyle \mathcal{P}_{bd}(E)& = &\big\{S \in \mathcal{P}(E):\, S \hbox{ is bounded}\big\},\\
\displaystyle \mathcal{P}_{cl,cv}(E) &=& \big\{S \in \mathcal{P}_{cl}(E):\, S \hbox{ is convex }\big\}.
%
\end{eqnarray*}
 Let $S$ be a nonempty subset of a Banach space $E$ and let $F : S \longrightarrow \mathcal{P}(E)$ be a multi-valued
mapping.
 %
  For every subset $M$ of $E$, we
write $$F^{-1}(M) = \left\{x \in S :\, F(x) \cap M\neq \emptyset \right\}$$ and for every subset $M$ of $S$ we write $$F(M)=\bigcup_{x\in M}F(x).$$
\begin{defin} Let $F : S \rightarrow \mathcal{P}(E)$ be a multi-valued operator. We say that:
\par \vskip0.1cm
\noindent $(i)$ $F $ has a weakly sequentially  closed graph if for every sequence $\{x_n\}_{n=1}^\infty$ of element in $S$ such that $x_n\rightharpoonup x$ in $S$ and for every sequence $\{y_n\}_{n=1}^\infty$ of elements in $E$ with $y_n\in F(x_n)$ such that $y_n\rightharpoonup y$ in $E$, then $y\in F(x)$. If $F$
is a single-valued mapping, then $F$ is called sequentially weakly continuous if for
any sequence $\{x_n\}_{n=1}^\infty$ in $S$ such that $x_n \rightharpoonup x \in S$, then $F(x_n) \rightharpoonup F(x).$
\par \vskip0.1cm
\noindent $(ii)$ $F$ is weakly
compact, if $F(A)$ is a relatively weakly compact, for all bounded subset $A$ of $S.$
\par \vskip0.1cm
\noindent $(iii)$
$F$ is sequentially weakly upper semi-compact in $S$ $($s.w.u.sco. for short$)$ if for any weakly convergent
sequence $\{x_n\}_{n=1}^\infty$ of elements in $S$ and for any arbitrary
 $y_n\in F(x_n)$, the sequence
 $\{y_n\}_{n=1}^\infty$ has a weakly convergent subsequence.
$\hfill\lozenge$\end{defin}

\bigskip

\noindent The theory of measures of weak noncompactness is an important tool used in
this work. This measure was introduced by De Blasi in
 \cite{De Blasi}. It is defined in the following way:
$$\beta(W):=\inf\big\{\varepsilon>0; \hbox{ there exist } K\in \mathcal W(E) \hbox{ and } \varepsilon>0 \hbox{ such that }W\subset K+\mathcal{B}_E(\varepsilon)\big\} $$
for all $W \in \mathcal{B}(E).$
\\
\noindent This measure share several properties such as
the subadditivity, the maximality, the positive homogeneousity, the monotonicity, and the nonsingularity (see \cite{De Blasi}).

\bigskip

Next, we state some definitions and a detailed information can
be found in \textup{\cite{amar2018,In}}.
\begin{defin}
 Let $F : \Omega \longrightarrow \mathcal{P}(X)$ be a multi-valued mapping. We say that:
  \par \vskip0.1cm
 \noindent $(i)$ $F$ is $\mathcal{D}$-set-Lipschitzian (with respect to $\beta$), if there is a continuous nondecreasing function $\phi : \mathbb{R}_+ \longrightarrow \mathbb{R}_+$ with $\varphi(0) = 0$ such that for any  bounded subset $W$ of $S$ with  $F(W)\in \mathcal{P}_{bd}(E),$ we have
 $$\beta(F(W ))\leq\phi(\beta(W )).$$ If $\phi(r)< r$ for $r>0$, then $F$ is called nonlinear $\mathcal{D}$-set-contraction.
 Moreover, if $\phi(r)=k\, r$ with $k<1,$ we say that $F$ is $k$-set-contraction.
 \par \vskip0.1cm
\noindent $(ii)$ $F$ is $\beta$-condensing, if $F$ is bounded and for any bounded subset $W$ of $\Omega$ with $\beta(W)> 0,$ we have $$\beta(F(W)) < \beta(D).$$
\par \vskip0.1cm
\noindent $(iii)$  $F$ is   $\beta$-nonexpansive, if $F$ is bounded and for any  bounded subset $W$ of $S,$ we have $$\beta(F(W)) \leq\beta(W).$$
~~$\hfill\lozenge$
 \end{defin}

 \bigskip

 \begin{defin}
 Let $F : S \longrightarrow \mathcal{P}(X)$ be a multi-valued mapping. We say that:
  \par \vskip0.1cm
 \noindent $(i)$ $F$ is countably $\mathcal{D}$-set-Lipschitzian,
  if there is a continuous nondecreasing function $\phi : \mathbb{R}_+ \longrightarrow \mathbb{R}_+$ with $\Phi(0) = 0$ such that for all countably bounded subset $W$ of $S$ with $F(W)\in \mathcal{P}_{bd}(E),$ we have
  $$\beta(F(W )) \leq\Phi(\beta(W )).$$
 In addition, if $\phi(r)< r,\,r>0,$ then
  $F$ is called countably $\mathcal{D}$-set-contraction.
 \par \vskip0.1cm
\noindent $(ii)$  $F$ is countably $k$-set-contraction, $0\leq k<1,$
 if for every countable bounded subset $W$ of $S$ with $F(W)\in \mathcal{P}_{bd}(E),$ we have $$\beta(F(W ))\leq k\,\beta(W )).$$
\par \vskip0.1cm
\noindent $(iii)$  $F$ is countably $\beta$-condensing,
if $F(S)$ is bounded and for every countable bounded subset $W$ of $S$ with $\beta(W)> 0,$ we have $$\beta(F(W)) < \beta(W).$$
\par \vskip0.1cm
\noindent $(iv)$  $F$ is countably  $\beta$-nonexpansive,
if $F$ is bounded and   for every countable bounded subset $W$ of $S,$ we have $$\beta(F(W)) \leq\beta(W).$$~~$\hfill\lozenge$
 \end{defin}

\bigskip

\noindent Because it lacks the stability of convergence for the product sequences under the
weak topology, A. Ben Amar, S. Chouayekh and A. Jeribi have introduced, in \cite{ASA}, a class of Banach
algebras satisfying a certain sequential condition $(\mathcal{P}):$
$$(\mathcal P)\left\{  \begin{array}{ll}
    \hbox{for any sequences } \{\alpha_n\}_{n=1}^\infty\hbox{ and }\{\beta_n\}_{n=1}^\infty\hbox{ of }X \hbox{ such that }\alpha_n\rightharpoonup \alpha
\hbox{ and } \beta_n \rightharpoonup \beta,\\ \hbox{ then } \alpha_n\cdot \beta_n \rightharpoonup \alpha\cdot \beta.
  \end{array}
\right.
$$
This class includes the finite dimensional Banach algebra and the space $C(K, X),$ where $K$ is a compact Hausdorff space and $X$ satisfying $(\mathcal{P})$.
 In order to prove some fixed point theorems
in Banach algebras satisfying $(\mathcal{P}),$ the authors in \cite{BA2016} have  introduced the concept of multi-valued mappings
of the form $\left(\frac{I-C}{A}\right),$ where $A$ and $C$ define multi-valued mappings acting
on Banach algebras.
\begin{defin}\label{reg}\cite{BA2016}
 Let $E$ be a Banach algebra and let $A, C : E \longrightarrow \mathcal{P}(E)$ be
multi-valued mappings. We say that the mapping $\frac{I-C}{A}$ is well defined on $x \in E$
and we write
$$y\in \left(\frac{I-C}{A}\right)(x)$$
if $x \in yA(x) + C(x).$ $\hfill\lozenge$
\end{defin}

\bigskip

\noindent Following \cite{Cardinali}, we recall the next notion.
\begin{defin}  Let $F : S \rightarrow \mathcal{P}(E)$ be a multi-valued operator. We say that $F $ has weakly closed graph in $S\times E,$ if for every sequence $\{x_n,n\in \mathbb{N}\}\subset S$ such that $x_n\rightarrow x,$ $x\in S$ and $y_n\in Fx_n$ such that $y_n \rightarrow y$ then $Fx\cap L(x,y)\neq \emptyset,$ where $L(x,y):=\{\alpha y+(1-\alpha) x, \alpha\in [0,1]\}.$\\
We say that $F $ has $w$-weakly closed graph in $S\times E,$
if it has weakly closed graph in $S\times E$ with respect to the weak topology.$\hfill\lozenge$
\end{defin}

\bigskip

\noindent Using the concept of multi-valued mappings with $w$-weakly closed graphs, A. Ben Amar et al. have established the following useful result.
\begin{thm}\cite{amar2018}\label{cor} Let $K$ be a nonempty, closed and convex subset of a Banach space $E$ and let $F : K \longrightarrow \mathcal{P}_{cl,cv}(K)$ be a multi-valued operator such that:
\par \vskip0.1cm
\noindent $(i)$ $F$ maps weakly compact sets into weakly relatively compact sets, and
\par \vskip0.1cm
\noindent $(ii)$ $F$ is countably $\omega$-condensing with $w$-weakly closed graph.
\par \vskip0.1cm
\noindent Then $F$ has a fixed point. $\hfill\lozenge$\end{thm}
\section{\textbf{Fixed point results for multi-valued mappings in Banach spaces}}
In this section, we establish some fixed point results for countably
$\mathcal{D}$-set-Lipschitzian and countably $\beta$-nonexpansive multi-valued mappings in Banach spaces under a weak topology setting. Let $S$ be a nonempty closed convex subset of a Banach space $E.$ To any multi-valued mapping $T : S\rightarrow \mathcal{P}(S)$ we associate a sequence $\left\{S^T_n,\, n\in\mathbb{N}\right\}$ of $\mathcal{P}_{cl,cv}(S)$
 defined by
\begin{eqnarray}\label{S}S^T_0= {S}\hbox{ and }
    S^T_{n+1}=\overline{co}\left\{T\left(S^T_n\right)\right\}, \hbox{ for } n=0, 1,\ldots.\end{eqnarray}
It is well known that $\left\{S^T_n,\, n\in\mathbb{N}\right\}$ is a decreasing sequence of nonempty subsets of $S.$
In the following result we give a fixed point theorem on some class of multi-valued mappings, as special cases, includes the multi-valued mappings with weakly sequentially closed graphs, which are countably $\beta$-condensing, countably $\mathcal{D}$-set-contraction, with respect to the De Blasi measure of weak noncompactness, and hemi-weakly compact, countably $\beta$-nonexpansive.
\begin{thm}\label{thm1} Let $S$ be a nonempty closed convex subset of $E$ and let $T : S\rightarrow \mathcal{P}(S)$ be a multi-valued mapping with a weakly sequentially closed graph  such that:
\par \vskip0.1cm
\noindent $(i)$ There exists a bounded subset $\Omega$ belongs to $\left\{S^T_n,\, n\in\mathbb{N}\right\},$
\par\vskip0.1cm
\noindent $(ii)$ $T$ is hemi-weakly compact and has closed, convex values on $\Omega,$
\par\vskip0.1cm
\noindent $(iii)$ $T$ is countably $\mathcal{D}$-set-Lipschitzian on $\Omega$ with $\mathcal{D}$-function $\varphi.$
\par\vskip0.1cm
\noindent Then $T$ has,  at least,  one fixed point in $\Omega$ as soon
as $\varphi(r) \leq r$ for $0< r\leq \|\Omega\|.$ $\hfill\lozenge$
\end{thm}
\noindent\textbf{Proof.} Let $n_0\in \mathbb{N}$ such that $\Omega=S^T_{n_0}.$
 Let $\gamma\in(0,1)$ and $y$ be fixed in $T^{n_0}(S).$ We define a multi-valued mapping  $T_{\gamma}: S\longrightarrow \mathcal{P}(S)$ by the formula $$T_{\gamma}(x)=\gamma T(x) +(1-\gamma) y.$$
Obviously $\left\{T^n(S),\, n\in\mathbb{N}\right\}$ is a decreasing sequence of nonempty subsets of $S$ since $S\neq \emptyset$ and $T(S)\subset S.$
In view of
the convexity of $S$, this leads to the conclusion that $ T_{\gamma}(S)\subset S.$   Consequently,  $\left\{S_n^{T_\gamma},\, n\in\mathbb{N}\right\}$ is, also, a decreasing sequence of nonempty subsets of $S.$ In particular, we have $$T_{\gamma}\left(S^{T_{\gamma}}_{n_0}\right)\subset
\overline{co}\left\{T_{\gamma}\left(S^{T_{\gamma}}_{n_0}\right)\right\}\subset S^{T_{\gamma}}_{n_0}.$$ Thus,
$T_{\gamma}$ maps  $S^{T_{\gamma}}_{n_0}$ into $\mathcal{P}\left(S^{T_{\gamma}}_{n_0}\right).$
Now, let us prove that $T_{\gamma}$ has convex values on $S^{T_{\gamma}}_{n_0}.$ To do this, we claim first that
 \begin{eqnarray}\label{3}S^{T_{\gamma}}_k\subset S^{T}_k \hbox{ for all }k\in \{0,1,\ldots, n_0\}.\end{eqnarray}
 We proceeding by induction argument. Evidently the inclusion (\ref{3}) hold for $K=0.$ Assume that  $S^{T_\gamma}_j\subset S^{T}_j$  for some $J\in \{1,\ldots, n_0-1\}.$
 Then,
 $$\begin{array}{rcl}\label{Qy}S^{T_\gamma}_{j+1}&= & \overline{co}\left\{T_\gamma\left(S^{T_\gamma}_{j}\right)\right\}\\
 &= & \overline{co}\left\{\gamma T\left(S^{T_\gamma}_{j}\right)+(1-\gamma)y\right\}\\
  &\subset & \overline{co}\left\{\gamma T\left(S^{T}_{j}\right)+(1-\gamma)y\right\}\\
    &\subset & \overline{co}\left\{\gamma\overline{co}\left(T\left(S^{T}_{j}\right)\right)
    +(1-\gamma)y\right\}\\
        &= & \overline{co}\left\{\gamma S^{T}_{j+1}+(1-\gamma)y\right\}.
\end{array}$$
Since $y\in T^{j+1}(S)\subset S^{T}_{j+1}$ and the set $S^{T}_{j+1}$ is convex and closed, we reach to result that $$\overline{co}\left\{\gamma S^{T}_{j+1}+(1-\gamma)y\right\}\subset \overline{co}\left\{\gamma S^{T}_{j+1}+(1-\gamma)S^{T}_{j+1}\right\}\subset S^{T}_{j+1}.$$ This proves the claim, and achieves that $S^{T_\gamma}_{n_0}$ is bounded. Now, let $x$ be an arbitrary element of $S^{T_{\gamma}}_{n_0}$ and let $t\in (0,1)$ and  $y_1, y_2\in T_\gamma(x).$ Then, there exist $u_1, u_2 \in T(x)$ such that $$y_1 =\gamma u_1+(1-\gamma)y\hbox{ and } y_2 =\gamma u_2+(1-\gamma)y.$$
Then, we have
$$t y_1+(1-t)y_2=  \gamma(tu_1+(1-t)u_2)+(1-\gamma)y.$$
Since $T$ has convex values, we deduce that $t y_1+(1-t)y_2\in T_\gamma(x).$ This prove the claim. Now let us show that $T_{\gamma}$ has closed values on $S^{T_{\gamma}}_{n_0}.$ Let $x\in S^{T_{\gamma}}_{n_0}$ and $\{z_n\}_{n=0}^\infty\subset T_{\gamma}(x)$ be a convergent sequence to some $z\in S^{T_{\gamma}}_{n_0}.$ Then there exist a sequence $\{w_n\}_{n=0}^\infty\subset T(x)$ such that $$z_n=\gamma w_n+(1-\gamma)y\hbox{ for all } n\in \mathbb{N}.$$
This implies that $w_n\rightarrow \frac{1}{\gamma}(z+(1-\gamma)y).$ Since $T$ has closed values, we get $\frac{1}{\gamma}(z+(1-\gamma)y)\in Tx,$ that is $z\in T_\gamma(x).$ This achieves the proof of our claim.
 Now we shows that $T_{\gamma}$ is countably $\beta$-condensing on $S^{T_\gamma}_{n_0}.$
 Take an arbitrary  countably subset $M$ of $S^{T_\gamma}_{n_0}.$ By using the boundedness of $T_\gamma(M)$ together with  the properties of the De
Blasi measure of weak noncompactness, we get
 \begin{eqnarray}\label{Q6}\beta(T_\gamma(M))
\leq \beta(\gamma T(M)+(1-\gamma)y)\leq \gamma\varphi(\beta(M)).\end{eqnarray}
If $\beta(M)>0$ we get $\beta(T_\gamma(M))
<\beta(M),$ which means that $T_\gamma$ is countably $\beta$-condensing since $T_\gamma$ is bounded on $\Omega.$
 Our next task is to show that $T_\gamma$ maps weakly compact sets into relatively weakly compact sets. To do so, take  a weakly compact set
  $M$ of $S^{T_\gamma}_{n_0}.$ Let $\{y_n\}_{n=1}^\infty$ be a sequence of $T_\gamma(M).$ Then there exists a sequence $\{x_n\}_{n=1}^\infty\subset M$
   such that $$y_n\in T_{\gamma}(x_n), \, n=0,1, \ldots.$$
 The use of $(\ref{Q6})$ allows us to obtain that \begin{eqnarray*}\label{Q1}\beta\left(\left\{y_n\right\}_{n=1}^\infty\right)
 \leq\beta\left(T_\gamma\left(\left\{x_n\right\}_{n=1}^\infty\right)\right)
\leq \gamma\varphi\left(\beta\left(\left\{x_n\right\}_{n=1}^\infty\right)\right)=0.\end{eqnarray*} Hence and in view of the Eberlein-\v{S}mulian's theorem, we infer that $T_\gamma(M)$ is relatively weakly compact. Since $T_\gamma$ has a weakly sequentially closed graph, then we may invoke Theorem \ref{cor} in order to conclude that  there exists $u_\gamma\in S_{n_0}^{T_\gamma}$ such that $$u_\gamma\in T_\gamma(u_\gamma).$$
 Now, let $\{\gamma_m\}_{m=1}^\infty$ be a sequence in $(0, 1)$ such that $\gamma_m\rightarrow 1.$ From the above discussion and by using inclusion (\ref{3}), there exists a sequences $\{u_m\}_{m=1}^\infty\subset \Omega$ such that
$$u_m\in  T_{\gamma_m }(u_m).$$
This leads to the conclusion that there exists $w_m\in T(u_m)$ such that
$$u_m=\gamma_m w_m+(1-\gamma_m)y.$$
On the other hand, it follows from inclusions $$T\left(\{u_m\}_{m=1}^\infty\right)\subset T\left(\Omega\right)\subset S^{T}_{n_0+1}\subset \Omega,$$ that $T(\{u_m\}_{m=1}^\infty)$ is bounded in $\Omega.$
Then, we result that $$u_m-w_m= (\gamma_m-1)w_m+(1-\gamma_m)y \rightharpoonup \theta.$$
Keeping in mind the hemi-weak compactness of $T,$ we deduce that
  $\{u_m\}_{m=1}^\infty$ has a weakly convergent subsequence to some $u.$
  Since
$T$ has a weakly sequentially closed graph, we reach the result that $u \in T(u).$
~~$\hfill\Box$

\bigskip

As an immediately consequence of Theorem \ref{thm1} we reach
the following fixed point theorem.
 \begin{cor} Let $S$ be a nonempty closed convex subset of $E$ and let $T : S\rightarrow \mathcal{P}(S)$ be a multi-valued mapping with a weakly sequentially closed graph  such that:
\par\vskip0.1cm
\noindent $(i)$ There exists a bounded subset $\Omega$ belongs to $\left\{S^T_n,\, n\in\mathbb{N}\right\},$
\par\vskip0.1cm
\noindent $(ii)$ $T$ is hemi-weakly compact and has closed, convex values on $\Omega,$
\par\vskip0.1cm
\noindent $(iii)$ $T$ is countably $\beta$-nonexpansive on $\Omega.$
\par\vskip0.1cm
\noindent Then $T$ has,  at least,  one fixed point in $\Omega.$ $\hfill\lozenge$
\end{cor}

\bigskip

The following results give a new fixed point theorem for countably
$\mathcal{D}$-set-contraction multi-valued mapping with $\mathcal{D}$-contraction function $\phi$ satisfying the estimates $\phi(r)<r$ only on a bounded real interval.
\begin{cor}\label{cor5} Let $S$ be a nonempty closed convex subset of $E$ and let $T : S\rightarrow \mathcal{P}(S)$ be a multi-valued mapping with a weakly sequentially closed graph  such that:
\par\vskip0.1cm
\noindent $(i)$ There exists a bounded subset $\Omega$ belongs to $\left\{S^T_n,\, n\in\mathbb{N}\right\},$
\par\vskip0.1cm
\noindent $(ii)$ $T$ has closed, convex values on $\Omega,$
\par\vskip0.1cm
\noindent $(iii)$ $T$ is countably $\mathcal{D}$-set-contraction on $\Omega,$ with $\mathcal{D}$-function $\varphi.$
\par\vskip0.1cm
\noindent Then $T$ has,  at least,  one fixed point in $\Omega.$ $\hfill\lozenge$
\end{cor}
\noindent\textbf{Proof.} Let $n_0\in \mathbb{N}$ such that $\Omega=S^T_{n_0}.$ Thanks to Theorem \ref{thm1}, we have only to show that $T : \Omega\longrightarrow \mathcal{P}\left(\Omega\right)$ is hemi-weakly compact. To do this, let $\{x_n\}_{n=1}^\infty$ be an arbitrary sequence of $\Omega$ and let $y_n\in T(x_n)$ such that $\{x_n-y_n\}_{n=1}^\infty$ has a weakly convergent subsequence.
 Furthermore, taking into account that $$\{x_n\}_{n=1}^\infty\subset \{x_n-y_n\}_{n=1}^\infty+T\left(\{x_n\}_{n=1}^\infty\right)$$
and bearing in mind the subadditivity of the measure of weak noncompactness of De Blasi, we arrive at the following conclusion:
 $$\begin{array}{rcl}\beta\left(\{x_n\}_{n=1}^\infty\right)&\leq & \beta\left(\{x_n-y_n\}_{n=1}^\infty\right)+\beta\left(T\left(\{x_n\}_{n=1}^\infty\right)\right)\\
 &\leq &\beta\left(T\left(\{x_n\}_{n=1}^\infty\right)\right).
\end{array}$$
By using a contradiction argument we can prove that $\{x_n\}_{n=1}^\infty$ is  relatively weakly compact. Consequently, $\{x_n\}_{n=1}^\infty$ has a weakly convergent subsequence, in view of the Eberlein-\v{S}mulian's theorem, which shows that $T$ is hemi-weakly compact.~~$\hfill\Box$

\bigskip

\begin{rem}
Corollary \ref{cor5} extends Theorem $3.1$ in \cite{BA2016}, and shows that the condition "$\mathcal{D}$-set-contraction" can be relaxed by assuming that $T$ is countably $\mathcal{D}$-set-contraction, similarly the condition "$T(S)$ is bounded" can be relaxed by assuming that $S^T_{n_0}$ is bounded for some $n_0\in \mathbb{N}.$~~$\hfill\lozenge$
\end{rem}

\bigskip

The following result concerning the countably $\beta$-condensing multi-valued mappings and represents an extension of Theorem \ref{cor}, since it does not require the boundness of $T(S).$
\begin{cor}
 Let $S$ be a nonempty closed convex subset of $E$ and let $T : S\rightarrow \mathcal{P}(S)$ be a multi-valued mapping with a weakly sequentially closed graph  such that:
\par\vskip0.1cm
\noindent $(i)$ There exists a bounded subset $\Omega$ belongs to $\left\{S^T_n,\, n\in\mathbb{N}\right\},$
\par\vskip0.1cm
\noindent $(ii)$ $T$ has closed, convex values on $\Omega,$
\par\vskip0.1cm
\noindent $(iii)$ $T$ is countably $\beta$-condensing and maps weakly compact sets into weakly relatively compact sets
\par\vskip0.1cm
\noindent Then $T$ has,  at least,  one fixed point in $\Omega.$ $\hfill\lozenge$
\end{cor}
\noindent\textbf{Proof.} Let $n_0\in \mathbb{N}$ such that $\Omega=S^T_{n_0}.$
Since $T$ is countably $\beta$-condensing on $\Omega,$ then
\begin{eqnarray}\label{7}\beta(T(M))\leq \beta(M) \hbox{ for all } M\subset \Omega \hbox{ such that }\beta(M)>0.
\end{eqnarray}
On the other hand, taking into account that $T$ maps weakly compact sets into weakly relatively compact sets, we obtain that $\beta(T(M))=0$ for all weakly compact subset $M$ of $\Omega.$ In particular, for every relatively weakly compact subset $M$ of  $\Omega,$ i.e $\beta(M)=0,$ we have
\begin{eqnarray}\label{8}\beta(T(M))\leq \beta(T(\overline{M}^w)) =0.
\end{eqnarray}
Using inequalities (\ref{7}) and (\ref{8}) we infer that $T$ is countably $\mathcal{D}$-set-Lipschitzian with $\mathcal{D}$-function $\varphi,$ given by $\varphi(r)=r$ for all $0\leq r\leq \|\Omega\|.$ Now, by proceeding essentially as in the proof of Corollary \ref{cor5} we can obtain that $T$ is hemi-weakly compact.
Now, we may invoke Theorem \ref{thm1} in order to conclude that $T$ has a fixed point $x$ in $\Omega.$~~$\hfill\Box$

\bigskip

 Next, we establish some new hybrid fixed point theorems involving the sum of two multi-valued mappings defined on a nonempty, closed, and convex subsets of Banach spaces.
\begin{thm}\label{TT6}
Let $S$ be a nonempty, closed, convex subset of a Banach
space $E,$ $A, B: S \longrightarrow \mathcal{P}_{cl,cv}(E)$ be two multi-valued mappings have weakly sequentially closed graphs such that:
\par\vskip0.1cm
\noindent $(i)$ $A$ is hemi-weakly compact and countably $\mathcal{D}$-set-lipschitzian with $\mathcal{D}$-function $\Phi$,
    \par\vskip0.1cm
\noindent $(ii)$ $B$ is weakly compact,
\par\vskip0.1cm
\noindent $(iii)$ $(A+B)(S)$ is bounded in $S .$
\par\vskip0.1cm
\noindent Then, the operator inclusion $(u\in Au+Bu)$ has,  at least,  one fixed point in $S$ as soon as $\Phi(r) \leq r$ for $r>0.$  $\hfill\lozenge$
\end{thm}
\noindent\textbf{Proof.} Since $S$ is a nonempty convex subset of $E$ and $(A+B)(S)\subset S,$ then we may define a multi-valued mapping $T$ by
$$\left\{
                           \begin{array}{ll}
                              T : S_{1}^{A+B}\longrightarrow  \mathcal{P}_{cl,cv}\left(S_{1}^{A+B}\right)\\
                             u\mapsto \displaystyle Au+ Bu
                           \end{array}
                         \right.
$$
Recall that $S_{1}^{A+B}=\overline{co}(A+B)(S)$ is the second terns of the sequence (\ref{S}) associate to $A+B.$
 Let $M$ be a countably bounded subset of $ S_{1}^{A+B}.$ Keeping in mind the relatively weak compactness of $B(M)$ and using the subadditivity of the De Blasi measure of weak
noncompactness we get
 $$\begin{array}{rcl}\label{Qy}\beta(T(M))&\leq& \beta (A(M))+\beta(B(M))\\
&\leq& \Phi(\beta(M)).
\end{array}$$
 Then, $T$ defines a countably $\mathcal{D}$-set-Lipschitzian multi-valued mapping on $S_{1}^{A+B}$ with $\mathcal{D}$-function $\Phi.$
Now, we claim that $T$ is hemi-weakly compact. To do this, let
$\left\{u_n\right\}_{n=1}^{\infty}$ be a sequence of elements in $S_{1}^{A+B}$ and let $v_n\in T(u_n)$ such that $$w_n:=u_n-v_n\rightharpoonup w\in S_{1}^{A+B}.$$ Since $v_n\in T(u_n)$ then there exist $a_n\in A(u_n)$ and $b_n\in B(u_n)$ such that
$$v_n=a_n+b_n,$$
 which implies that, $$u_n-a_n=w_n+b_n.$$
It follows from the Eberlein-\v{S}mulian's theorem, that  $\{b_n,\ n\in \mathbb{N}\}$ has a weakly convergent subsequence. This shows from the hemi-weak compactness of $A,$ that $\{u_n,\ n\in \mathbb{N}\}$ has a weakly convergent subsequence.
It is easily seen that, $T$ has a weakly sequentially closed graph.
 So, $T$ has a fixed point in $S_{1}^{A+B}$
due to Theorem \ref{thm1}.
 $\hfill\Box$

\bigskip

Now, we are in a position to apply Theorem \ref{TT6} in order to
deduce the following fixed point theorem.
\begin{cor}\label{cor1}
Let $S$ be a nonempty, closed, convex subset of a Banach
space $E,$ $A, B: S \longrightarrow \mathcal{P}_{cl,cv}(E)$ be two multi-valued mappings
 have weakly sequentially closed graphs such that:
 \par\vskip0.1cm
 \noindent $(i)$ $A$ is hemi-weakly compact and countably $\beta$-nonexpansive,
\par\vskip0.1cm
\noindent $(ii)$ $B$ is weakly compact,
\par\vskip0.1cm
\noindent $(iii)$ $(A+B)(S)$ is a bounded subset of $S.$
\par\vskip0.1cm
\noindent Then, the operator inclusion $(u\in Au+Bu)$ has,  at least,  one fixed point in $S.$ $\hfill\lozenge$
\end{cor}

\bigskip

\begin{rem}
Theorem \ref{TT6} extends and improves Theorem $4.10$ in \cite{BADO2016}. Indeed, if there exist a bounded subset $S_0$ in $E$ and a sequence $\{\lambda_n,\ n\in \mathbb{N}\}\subset (0,1)$ with $\lambda_n \rightarrow 1$ such that $(A+\lambda_n B)(S)\subset S_0$ for all $n,$ then $S^{A+B}_1$ is bounded. $\hfill\lozenge$
\end{rem}
\section{\textbf{Fixed Point results for multi-valued mappings in Banach algebras}}
In this section, we give sufficient conditions for the operator inclusion (\ref{in0}) which is acting on a Banach algebra satisfying $(\mathcal{P})$ to have a fixed point. Firstly, we establish some hybrid fixed point
theorems for $(\ref{in0})$ by looking at the multi-valued mapping $$T:=A\cdot B+C.$$
\begin{thm}\label{TT5}
Let $S$ be a nonempty, closed and convex subset of a Banach
algebra $E$ satisfying $(\mathcal{P}).$ Suppose that $A, B, C : S \longrightarrow \mathcal{P}(E)$ are three multi-valued mappings
 have weakly sequentially closed graphs such that:
 \par\vskip0.1cm
 \noindent $(i)$ For each $u\in S,$ $A(u)\cdot B(u)+C(u)$ is a closed, convex subset of  $ S.$
 \par\vskip0.1cm
 \noindent $(ii)$ There exists a bounded subset $\Omega$ belongs to $\left\{S^{A\cdot B+C}_n,\, n\in\mathbb{N}\right\},$
 \par\vskip0.1cm
\noindent $(iii)$ $A, B$ and $C$ are countably $\mathcal{D}$-set-lipschitzian on $\Omega$ with $\mathcal{D}$-functions $\Phi_A, \Phi_B$ and $\Phi_C$ respectively,
\par\vskip0.1cm
\noindent $(iv)$ $A(\Omega)$ and $B(\Omega)$ are bounded.
\par\vskip0.1cm
\noindent Then, the operator inclusion $(\ref{in0})$ has,  at least,  one fixed point in $\Omega$ provided that\\ $\|B(\Omega)\|\Phi_A(r)+\|A(\Omega)\|\Phi_B(r)+\Phi_A(r)
\Phi_B(r)+\Phi_C(r)< r\hbox{ for }0< r\leq \|\Omega\|.$ $\hfill\lozenge$
\end{thm}
\noindent\textbf{Proof.} Let $n_0\in \mathbb{N}$ such that $\Omega= S^{A\cdot B+C}_{n_0}.$ Proceeding as in the proof of Theorem \ref{TT6}, in view of assumption $(i),$ we can define a multi-valued mapping $T$ by
$$\left\{
                           \begin{array}{ll}
                              T : \Omega\longrightarrow  \mathcal{P}_{cl,cv}(\Omega)\\
                             u\mapsto \displaystyle Au\cdot Bu+Cu.
                           \end{array}
                         \right.
$$
First, we claim that $T$ is countably $\mathcal{D}$-set-lipschitzian. Indeed,
 let $M$ be a countably bounded subset of $\Omega.$ Combining Theorem $2.9$ \cite{equivalence} with Lemma $
 3.1$ in \cite{Mefteh} we get $$\begin{array}{rcl}  \beta (A(M)B(M))
\leq \|B(M)\|\,\Phi_A(\beta(M))+\|A(M)\|\,\Phi_B(\beta(M))+\Phi_A(\beta(M))
\,\Phi_B(\beta(M)).
\end{array}$$
Then, keeping in mind the subadditivity property of $\beta,$ we have that
 $$\begin{array}{rcl}\label{Qy}\beta(T(M))&\leq&  \beta (A(M)B(M))+\beta(C(M))\\
&\leq& \|B(M)\|\,\Phi_A(\beta(M))+\|A(M)\|\,\Phi_B(\beta(M))+\Phi_A(\beta(M))
\,\Phi_B(\beta(M))+\Phi_C(\beta(M)).
\end{array}$$
 Thus, $T$ defines a countably $\mathcal{D}$-set-Lipschitzian on $\Omega$ with $\mathcal{D}$-function $\Psi$ given by $$\Psi(r)=\|B(\Omega)\|\,\Phi_A(r)+\|A(\Omega)\|\,\Phi_B(r)+\Phi_A(r)
\,\Phi_B(r)+\Phi_C(r)\hbox{ for }r\geq0.$$
Now, we claim that $T$ has a weakly sequentially closed graph. To this end, let  $\left\{u_n\right\}_{n=1}^{\infty}$ be a sequence of elements in $\Omega$ such that $u_n \rightharpoonup u\in \Omega$ and let $v_n\in T(u_n)$ such that $v_n\rightharpoonup v\in \Omega.$ Then there exist $\sigma_n\in Au_n$ and $\gamma_n\in Bu_n$ and $\rho_n\in Cu_n$ such that $$v_n=\sigma_n\cdot \gamma_n+\rho_n.$$
Since $A$ is countably $\mathcal{D}$-set-Lipschitzian, then
$$\beta\left(\left\{\sigma_n\right\}_{n=1}^{\infty}\right)\leq \Phi_A\left(\beta\left(\left\{u_n\right\}_{n=1}^{\infty}\right)\right)=0.$$
From the Eberlein-\v{S}mulian's theorem,
we can assume that $\left\{\sigma_n\right\}_{n=1}^{\infty}$ and $\left\{\gamma_n\right\}_{n=1}^{\infty}$ converge, respectively,  to some points $ \sigma, \gamma\in E.$ Further, taking into account that $A$ and $B$ have  weakly sequentially closed graphs  we get
 $\sigma\in Au$ and $\gamma\in Bu.$ From the condition $(\mathcal{P}),$ it follows that
 $$\rho_n \rightharpoonup v-\sigma\cdot \gamma.$$
Keeping in mind that $C$ has a weakly sequentially closed graph, we derive that $$v-\sigma\cdot \gamma\in Cx.$$
 Consequently, we get $$v\in Au\cdot Bu+Cu,$$ which achieves the proof of our claim.
 Now, an application of Corollary \ref{cor5} yields that there is a element $u\in \Omega$ such that $u\in Au\cdot Bu+Cu.$
 $\hfill\Box$

\bigskip

\begin{rem}
 Theorem \ref{TT5} extends Theorem $4.2$ in \cite{BA2016}, and shows that the condition "$\mathcal{D}$-set-Lipschitzian" can be relaxed by assuming that $F$ is countably $\mathcal{D}$-set-Lipschitzian, similarly the condition "$A(S), B(S)$ and $C(S)$ are bounded" can be relaxed by assuming that $S^{A\cdot B+C}_n$ is bounded for some $n\in \mathbb{N},$ moreover
it proves that the condition that $A, B$ and $C$ are s.w.u.sco., i.e. Condition $(3)$ in \cite[Theorem 4.2]{BA2016}, in the statement of that theorem is not needed. $\hfill\lozenge$
\end{rem}

\bigskip

\begin{rem}
 Since every s.w.u.sco., $k$-lipschitzian multi-valued mapping is $\mathcal{D}$-set-Lipschitzian, with $\mathcal{D}$-function $\phi(r)=k\,r,$ in particular is countably $\mathcal{D}$-set-Lipschitzian.  Moreover, by using the condition $(\mathcal{P}),$ it is clear that $A\cdot B+C$ has closed values if the multi-valued mappings $A, B,$ and $C$ have closed values. Thus Theorem \ref{TT5} generalizes Theorem $4.5$ in \cite{BA2016}, by taking $n_0=0.$
 $\hfill\lozenge$
\end{rem}

\bigskip

A special case of Theorem \ref{TT5}, which is useful in applications to differential and integral inclusions, is introduced in the following theorem.
\begin{cor}\label{cor6}
Let $S$ be a nonempty, closed and convex subset of a Banach
algebra $E$ satisfying $(\mathcal{P}).$ Suppose that $A, B, C : S \longrightarrow \mathcal{P}(E)$ are three multi-valued mappings
 have weakly sequentially closed graphs such that:
\par\vskip0.1cm
\noindent $(i)$ For each $u\in S,$ $A(u)\cdot B(u)+C(u)$ is a closed, convex subset of  $ S,$
 \par\vskip0.1cm
 \noindent $(ii)$ There exists a bounded subset $\Omega$ belongs to $\left\{S^{A\cdot B+C}_n,\, n\in\mathbb{N}\right\},$
 \par\vskip0.1cm
\noindent $(iii)$ $A$ and $C$ are countably $\mathcal{D}$-set-lipschitzian with $\mathcal{D}$-functions $\Phi_A$ and $\Phi_C$ respectively,
    \par\vskip0.1cm
\noindent $(iv)$ $B$ is weakly compact on $\Omega,$
\par\vskip0.1cm
\noindent Then, the operator inclusion $(\ref{in0})$ has,  at least,  one fixed point in $S$ provided that\\ $\|B(\Omega)\|\Phi_A(r)+\Phi_C(r)< r$ $\hbox{ for }0< r\leq \|\Omega\|.$ $\hfill\lozenge$
\end{cor}

\bigskip

\begin{rem}
 Corollary \ref{cor6} extends Theorem $4.3$ in \cite{BA2016}, and shows that the condition that $A$ and $C$ are s.w.u.sco., i.e. Condition $(4)$ in \cite[Theorem 4.3]{BA2016}, in the statement of that theorem is not needed.~~$\hfill\lozenge$
\end{rem}

\bigskip

When $A$ is countably $\beta$-nonexpansive we get the following result.
\begin{cor}\label{T2}
Let $S$ be a nonempty, closed and convex subset of a Banach
algebra $E$  satisfying $(\mathcal{P}).$ Suppose that $A, B, C : S \longrightarrow \mathcal{P}(E)$ are three multi-valued mappings have weakly sequentially closed graphs such that:
\par\vskip0.1cm
\noindent $(i)$ for each $u\in S,$ $A(u)\cdot B(u)+C(u)$ is a closed, convex subset of  $ S,$
\par\vskip0.1cm
 \noindent $(ii)$ There exists a bounded subset $\Omega$ belongs to $\left\{S^{A\cdot B+C}_n,\, n\in\mathbb{N}\right\},$
 \par\vskip0.1cm
\noindent$(iii)$ $A$ is countably $\beta$-nonexpansive and $B$ is weakly compact  on $\Omega,$
\par\vskip0.1cm
\noindent$(iv)$ $C$ is countably $\alpha$-set-contraction  on $\Omega,$
\par\vskip0.1cm
\noindent Then, the operator inclusion $(\ref{in0})$ has,  at least,  one fixed point in $S$ provided that  $\|B(\Omega)\|+\alpha<1.$ $\hfill\lozenge$
\end{cor}
\bigskip

The  operator $\left(\frac{I-C}{A}\right)$ was used in numerous works to establish  some existence results for Hybrid fixed point theorems for  several class of single-valued mappings in Banach algebras (see for example \textup{\cite{BA2010, ASA, BA regular}}).
In most of these works, the operator $\left(\frac{I-C}{A}\right)$ play a fundamental role in their arguments.
Recently, A. Ben amar et al have introduced in \cite{BA2016} an analogues definition of $\left(\frac{I-C}{A}\right)$ for the multi-valued mappings, in order to investigate some fixed point theorems  for (\ref{in0}).
Next, we consider Hybrid fixed point
theorems for $(\ref{in0})$ by looking at the multi-valued mapping $$T:=\left(\frac{I-C}{A}\right)^{-1}B$$
where $A, B$ and $C$ are countably
$\mathcal{D}$-set-Lipschitzian or countably nonexpansive multi-valued mappings.
\begin{thm}\label{TT1}
Let $S$ be a nonempty, closed, and convex subset of a Banach algebra
$E$ satisfying $(\mathcal{P}).$ Assume that $A, C : E \longrightarrow \mathcal{P}_{cl,cv}(E)$ and $B : S \longrightarrow \mathcal{P}_{cl,cv}(E)$ are three multi-valued mappings have weakly sequentially closed graphs such that:
\par\vskip0.1cm
\noindent $(i)$  $A$ is  weakly compact,
\par\vskip0.1cm
\noindent $(ii)$   $B$ is countably $\mathcal{D}$-set-contraction with $\mathcal{D}$-contraction $\Phi$,
    \par\vskip0.1cm
 \noindent $(iii)$   $C$ is $\mathcal{D}$-set-contraction with $\mathcal{D}$-function $\Psi$,
     \par\vskip0.1cm
 \noindent $(iv)$   $A(E), B(S)$, and $C(E)$ are bounded,
 \par\vskip0.1cm
\noindent $(v)$ for each $u\in S,$ $\left(\frac{I-C}{A}\right)^{-1}B(u)$ is a closed, convex subset of $S.$
\par\vskip0.1cm
  \noindent Then, the operator inclusion $(\ref{in0})$ has a solution in $S$ whenever $\|A(S)\|\,\Phi(r)+\Psi(r)< r\hbox{ for }r>0.$ $\hfill\lozenge$
  \end{thm}
\noindent\textbf{Proof.}
Let $\gamma\in B(S)$ be fixed. From our assumptions we can define a multi-valued mapping $Q$ by
$$\left\{
                           \begin{array}{ll}
                              Q: E^{A\cdot \gamma+C}_1 \longrightarrow \mathcal{P}_{cl,cv}\left(E^{A\cdot \gamma+C}_1\right)\\
                             x\mapsto Ax\cdot \gamma+Cx.
                           \end{array}
                         \right.
$$
Recall that the bounded subset $E^{A\cdot \gamma+C}_1=\overline{co}(A\cdot \gamma+C)(E)$ is the second terns of the sequence (\ref{S}) associate to $A\cdot \gamma+C.$
Proceeding essentially as in the proof of Theorem \ref{TT5}, we can show that $Q$ defines a countably $\mathcal{D}$-set-Lipschitzian, with $\mathcal{D}$-function $\Psi.$
Now we claim that $Q$ has a weakly sequentially closed graph. To do this, let $\{u_n\}_{n=1}^\infty$ be a sequence of elements in $E^{A\cdot \gamma+C}_1$ such that $u_n\rightharpoonup u$ and let $v_n\in Q(u_n)$ such that $v_n\rightharpoonup v.$
 Then there exist two sequences $\{\sigma_n\}_{n=1}^\infty$ and  $\{\rho_n\}_{n=1}^\infty$ with $\sigma_n\in Au_n$ and $\rho_n\in Cu_n$ such that $$v_n =\sigma_n\cdot \gamma+\rho_n.$$
 Since $\beta\left(A\left(\{u_n\}_{n=1}^\infty\right)\right) =0,$
taking into account the Eberlein-\v{S}mulian's theorem we can suppose that $\sigma_n \rightharpoonup \sigma\in E.$
 Keeping in mind that  $A$ and $C$ have weakly sequentially closed  graphs, we deduce that $\sigma\in Au,$ and  $ v-\sigma\cdot \gamma\in Cu.$
 This means that $v\in Tu.$  Invoking Corollary \ref{cor5}, we infer that $Q$ has a fixed point $u\in E^{A\cdot \gamma+C}_1,$ i.e. $$u\in Au\cdot \gamma+Cu.$$ Thus, $$\gamma\in \left(\frac{I-C}{A}\right)(u),$$ and so $$\left(\frac{I-C}{A}\right)(u)\cap B(S)\neq \emptyset.$$ This achieves that $\left(\frac{I-C}{A}\right)^{-1}$ is well defined on $B(S).$ Accordingly, we can define a nonempty closed convex subset $S^{T}_1:=\overline{co}\left(\left(\frac{I-C}{A}\right)^{-1}B(S)\right)$ of $S,$ which is invariant by $T:=\left(\frac{I-C}{A}\right)^{-1}B$ in view of assumption $(v).$
Using assumption $(iv)$ together with the inclusion  \begin{eqnarray}\label{Eq3}
T(D)\subset A T(D)\cdot B(D)+ C T(D) \hbox{ for all }D\subset S,
  \end{eqnarray}
  in order to results that $S^{T}_1$ is a bounded subset of $S.$
 Now, let us consider the following multi-valued mapping:
$$\left\{
                           \begin{array}{ll}
                             T: S^{T}_1\longrightarrow \mathcal{P}_{cl,cv}(S^{T}_1)\\
                             u\mapsto \displaystyle \left(\frac{I-C}{A}\right)^{-1}B(u)
                           \end{array}
                         \right.
$$
  We will prove that $T$ is countably condensing. To see this,
  let $D$ be a countably subset of $S^{T}_1.$ By combining Theorem $2.9$ \cite{equivalence} with Lemma $
 3.1$ in \cite{Mefteh}, and using the subadditivity of $\beta$ and the relatively weak compactness of  $A(S^{T}_1)$ it follows that \begin{eqnarray}\label{eq4}\begin{array}{rcl}\beta(T(D))& \leq& \|A(S^{T}_1)\|\,\beta (B(D))+\beta (C(T(D)))\\
 & \leq& \|A(S)\|\,\Phi\left(\beta(D)\right)+\Psi\left(\beta(T(D))\right).\end{array} \end{eqnarray}
 If $\Phi(\beta(D))=0,$ by using a contradiction argument we can show that $T(D)$ is relatively weakly compact.
So, we may assume that $\Phi(\beta(D))>0.$ Inequalities (\ref{eq4}) imply that
$$\beta(T(D))-\Psi\left(\beta(T(D))\right) \leq \|A(S)\|\,\Phi\left(\beta(D)\right).$$
Combining this inequality with our assumptions, we get
 $$\Phi(\beta(T(D)))<\Phi(\beta(D)).$$
This means that $T$ is countably condensing. In view of
inequalities (\ref{eq4}) and by using a contradiction argument, this leads to the conclusion that
  $T$ maps weakly compact sets into relatively weakly compact sets.
 Now we claim that $T$ has a weakly sequentially closed graph. To do this, let $\{u_n\}_{n=1}^\infty$ be a sequence of elements in $S^{T}_1$ which converges weakly to $ u$ and $v_n\in T(u_n)$ which converges weakly to $v.$ Then, there exists a sequence $\{w_n\}_{n=1}^\infty$ in $E$ such that
 \begin{eqnarray}\label{incl}w_n\in \left(\frac{I-C}{A}\right)(v_n)\cap B(u_n), \hbox{ for all }n\in \mathbb{N}.\end{eqnarray}
  Using the fact that $\{u_n\}_{n=1}^\infty$ is weakly convergent we get
$$\beta\left(\{w_n\}_{n=1}^\infty\right)\leq\beta\left(B\left(\{u_n\}_{n=1}^\infty\right)\right)
  \leq \Phi\left(\beta\left(\{u_n\}_{n=1}^\infty\right)\right)=0.$$
 Then, we can extracts a renamed subsequence of $\{w_n\}_{n=1}^\infty $  such that $w_n\rightharpoonup w,$ which is
a consequence of the Eberlein-\v{S}mulian's theorem. Keeping in mind that $B$  has a weakly sequentially closed graph, we derive that
 $$w\in B(u).$$ On the other hand, by using the Definition \ref{reg} together with inclusion (\ref{incl}), we deduce that $$v_n\in w_n \cdot A(v_n)+C(v_n), n\in \mathbb{N}.$$ Since $A$ and $C$ have weakly closed graphs and $E$ satisfying $(\mathcal{P}),$ we infer that $$v\in w \cdot A(v)+C(v),$$ and consequently $$w\in \left(\frac{I-C}{A}\right)(v).$$
 This means that $v\in T(u),$ and the claim is approved.
 The remained proof follows along the lines of Theorem \ref{cor}.
 \hfill $\Box$\par\medskip

\bigskip

\begin{thm}\label{TT4.5}
Let $S$ be a nonempty, closed and convex subset of a Banach
algebra $E$ satisfying $(\mathcal{P}).$ Assume that $A, C : E \longrightarrow \mathcal{P}_{cl,cv}(E)$ and $B : S \longrightarrow \mathcal{P}(E)$ are three multi-valued mappings have weakly sequentially closed graphs such that:
 \par\vskip0.1cm
\noindent$(i)$   $A$ is $\mathcal{D}$-set-Lipschitzian with $\mathcal{D}$-function $\Phi$,
    \par\vskip0.1cm
\noindent$(ii)$  $B$ is countably $\alpha$-set-contraction,
\par\vskip0.1cm
\noindent$(iii)$   $C$  is $\gamma$-set-contraction,
\par\vskip0.1cm
\noindent$(iv)$  $A(E),$ $B(S)$ and $C(E)$ are bounded,
\par\vskip0.1cm
\noindent$(v)$  for each $u\in S,$ $\left(\frac{I-C}{A}\right)^{-1}B(u)$ is a closed, convex subset of $S.$
\par\vskip0.1cm
\noindent Then, the operator inclusion $(\ref{in0})$ has a fixed point in $S$ whenever there exists a positif constant $\delta, \alpha<\delta<\frac{1-\gamma}{\|A(S)\|},$ such that $\displaystyle \|B(S)\|\Phi(r)<(1-2\delta\|A(S)\|-\gamma)r$ $ \hbox{ for } r>0.$
$\hfill\lozenge$
\end{thm}
\noindent\textbf{Proof.}
We shows, first, that
$\left(\frac{I-C}{A}\right)^{-1}$ exists on $B(S).$ To do this,
 let $v\in B(S)$ be fixed. Define a multi-valued mapping $Q$ by
$$\left\{
                           \begin{array}{ll}
                              Q : E \longrightarrow \mathcal{P}_{cl,cv}(E)\\
                             x\mapsto Ax\cdot v+Cx.
                           \end{array}
                         \right.
$$
An argument similar to that in the proof of Theorem \ref{TT1} yields that $Q$ is countably $\mathcal{D}$-set-contraction, with $\mathcal{D}$-function $\varphi(r)=\|v\|\Phi(r)+\gamma r,$
  and has a weakly sequentially closed graph.
On the other hand, by using a contradiction argument, we may prove that $Q$ maps weakly compact sets into weakly relatively compact sets.
 Applying now Theorem \ref{cor} we obtain that there exists $u\in E$ such that $$u\in Au\cdot v+Cu,$$ thus $$v\in \left(\frac{I-C}{A}\right)(u).$$
 As a result,
  $$\left(\frac{I-C}{A}\right)(u)\cap B(S)\neq \emptyset.$$
   This achieves that $\left(\frac{I-C}{A}\right)^{-1}$ is well defined on $B(S).$
Arguing as in the proof of Theorem \ref{TT1}, we can define
  a multi-valued mapping $T$ by
$$\left\{
                           \begin{array}{ll}
                              T : S_1^T\longrightarrow  \mathcal{P}_{cl,cv}\left(S_1^T\right)\\
                             u\mapsto \displaystyle \left(\frac{I-C}{A}\right)^{-1}B(u)
                           \end{array}
                         \right.
$$
 Proceeding essentially as in the proof of Theorem \ref{TT1}, we can prove that $T$ has a weakly sequentially closed graph.
Our next task is to show that $T$ is countably $\frac{\alpha}{\delta}$-set-contraction. To achieves this, take an arbitrary countably bounded subset $D\subset S_1^T$. Keeping in mind the inclusion (\ref{Eq3}) together with assumption $(iv)$ we infer that $S_1^T$ is a bounded subset of $S.$ Using again inclusion (\ref{Eq3}) combining with Theorem $2.9$ \cite{equivalence} with Lemma $
 3.1$ in \cite{Mefteh} we obtain
 \begin{eqnarray}\begin{array}{ll}\label{eq5}\beta(T(D))
&\leq \|B(S)\|\,\Phi(\beta(T(D)))+\alpha\,\|A(S)\|\, \beta(D)+\alpha\, \, \beta(D)\,\beta(A(T(D)))+\gamma\,\beta(T(D))\\
&\leq \|B(S)\|\,\Phi(\beta(T(D)))+2\alpha\,\|A(S)\|\, \beta(D)+\gamma\,\beta(T(D)),\end{array}\end{eqnarray}
 which implies that
 \begin{eqnarray*}(1-\gamma)\,\beta(T(D))
&\leq& \|B(S)\|\,\Phi(\beta(T(D)))+2\alpha \,\|A(S)\|\,\beta(D).\end{eqnarray*}
Then, by using the inequality $$\|B(S)\|\Phi(r)<(1-(2\,\delta\,\|A(S)\|+\gamma))\,r\hbox{ for }r>0,$$  we get $$\beta(T(D))\leq\frac{\alpha}{\delta}\beta(D).$$ This implies that $T$ is countably $\frac{\alpha}{\delta}$-set-contraction.
By proceeding essentially as in the proof of Theorem \ref{TT1}, we can prove that $T$ has a weakly sequentially closed graph.
 Hence, the desired result follows from Corollary \ref{cor5}. $\hfill\Box$

\bigskip

\begin{cor}
Let $S$ be a nonempty, closed and convex subset of a Banach
algebra $E$ satisfying the sequential condition $(\mathcal{P}).$ Assume that $A, C : E \longrightarrow \mathcal{P}_{cl,cv}(E)$, $B : S \longrightarrow \mathcal{P}(E)$ are three multi-valued mappings have weakly sequentially closed graphs such that:
\par\vskip0.1cm
\noindent $(i)$ $A$ is countably $\beta$-nonexpansive,
\par\vskip0.1cm
\noindent $(ii)$ $B$ is countably $\alpha$-set-contraction,
\par\vskip0.1cm
\noindent $(iii)$ $C$ is $\gamma$-set-contraction,
\par\vskip0.1cm
\noindent $(iv)$ $A(E),$ $B(S)$ and $C(E)$ are bounded,
\par\vskip0.1cm
\noindent $(v)$  for each $u\in S,$ $\left(\frac{I-C}{A}\right)^{-1}B(u)$ is a closed, convex subset of $S.$
\par\vskip0.1cm
\noindent Then, the operator inclusion $(\ref{in0})$ has a fixed point in $S$ whenever
 whenever there exists a positif constant $\delta, \alpha<\delta<\frac{1-\gamma}{\|A(S)\|},$ such that $2\delta\|A(S)\|+ \|B(S)\|+\gamma<1.$
$\hfill\lozenge$
\end{cor}
\section{\textbf{Nonlinear integral inclusions}}
Let $(X,\|\cdot\|)$ be a Banach algebra, we denote by $E:= C(J, X)$ the Banach algebra of all continuous functions from $J$ to $X,$ endowed with the supremum norm $$\|f\|_{\infty}=\sup\left\{\|f(t)\|~,~t\in J\right\}, \hbox{ for each } f\in E.$$
\noindent For any multi-valued mapping $\Xi : J\times X\longrightarrow \mathcal{P}(X)$ and for any $x \in E$ we denote
$$\displaystyle\int_0^t \Xi(s,x(s))\, ds = \left\{\int_0^t v(s)\, ds:\, v \hbox{ is Pettis integrable and } v(t)\in \Xi(t,x(t))\right\}$$
and
$$\|\Xi(t,x)\|=\sup\big\{\|v\| :\, v\in \Xi(t,x)\big\}.$$

\noindent In this section, we are mainly concerned with the existence results of solutions
for the following problem of nonlinear integral inclusion on $E$:
\begin{eqnarray}\label{NII}
   \displaystyle x(t)\in \displaystyle  T_1(t,x(t))  \cdot\big(q(t)+\int_0^{t}k(t,s)\,H(s,x(s))\, ds\big)+T_2(t,x(t)), \, t\in J.
\end{eqnarray}
where $J=[0,1],$ $k : J\times J \longrightarrow \mathbb{R}$, $q :J \longrightarrow X,$ $T_i : J \times X \longrightarrow X$ and $ H : J \times X \longrightarrow \mathcal{P}_{cp,cv}(X).$

\noindent Under the following assumptions, we could reach the solution of $(\ref{NII}):$
\par \vskip0.1cm
\begin{itemize}
\item[$(\bf{C_0})$] The functions $q: J\longrightarrow X$ and $ k: J\times J\longrightarrow \mathbb{R}$  are continuous.
    \par \vskip0.1cm
\item[$(\bf{C_1})$] For $i=1,2,$ we have:
\par \vskip0.1cm
\begin{itemize}
\item[$(\bf{i})$] For all $t\in J,$  $T_i(t,\cdot)$ are weakly sequentially  continuous.
    \par \vskip0.1cm
\item[$(\bf{ii})$] There exist two continuous functions $\vartheta_i : J \longrightarrow [0,\infty)$ such that
$$\big\|T_1(t, x) - T_1(t', y)\big\| \leq \vartheta_1(t)\,\|x - y\|+|\vartheta_1(t)-\vartheta_1(t')|\,\|y\|$$ $\hbox{ for all }x, y \in X \hbox{ and }t \in J,$
and $$\big\|T_2(t, x(t)) - T_2(t', x(t'))\big\| \leq |\vartheta_2(t)-\vartheta_2(t')| \, \|x\|_\infty$$ $\hbox{ for all }x \in E \hbox{ such that } \|x\|_\infty\leq R $ and $t, t'\in J.$
\par \vskip0.1cm
\item[$(\bf{iii})$] There exists a non decreasing mapping $\Psi : [0,+\infty) \longrightarrow [0,+\infty)$ such that  $$\beta(T_2(J\times M)) \leq \Psi(\beta (M))$$ $\hbox{ for all countably bounded subset }M \hbox{ of } X.$
\end{itemize}
\par \vskip0.1cm
\item[$(\bf{C_2})$] The multi-valued mapping $H : J\times X\longrightarrow \mathcal{P}(X)$ is such that:
    \par \vskip0.1cm
\begin{itemize}
\item[$(\bf{i})$] For each $x\in {E},$ there are a scalarly measurable function  $w : J \longrightarrow X$ with $w(t) \in H(t, x(t))$ a.e. $t\in J$ such that $w$ is Pettis integrable on $J$.
    \par \vskip0.1cm
\item[$(\bf{ii})$] There exists $R>0$ and   $h\in L^1(J,\mathbb{R}_+)$  such that $$ \|H(s, x)\|\leq h(s)\|x\|, \hbox{ for all } x\in \mathcal{B}_X(R). $$
    \par \vskip0.1cm
\item[$(\bf{iii})$]  There exists a non decreasing mapping $\Phi : [0,+\infty) \longrightarrow [0,+\infty)$ such that  $$\beta(H(J\times M(J)))  \leq \Phi(\beta (M(J)))$$ $\hbox{ for all countably bounded subset }M \hbox{ of } \mathcal{B}_E(R).$
    \par \vskip0.1cm
\item[$(\bf{iv})$] $H(t,\cdot)$ has weakly sequentially closed graph  on $\mathcal{B}_X(R).$
\end{itemize}
\end{itemize}

\bigskip

\noindent Before reaching the main result in this section, the useful results for the sequel are stated.
\begin{thm}\label{Mitchell}\cite{Mitchell}
Let $X$ be a Banach algebra and let $H\subset C([0,T],X)$ be bounded and equi-continuous. Then the map $t \longrightarrow \beta(H(t))$ is continuous on $[0,T]$ and
$$\beta(H) = \sup_{t\in[0,T]}\beta(H(t))=\beta(H([0,T])),$$
 where $H(t)=\{x(t)~,~x\in H\}$ and $H([0,T])=\bigcup_{t\in [0,T]}\{x(t)~,~x\in H\}.$ $\hfill\lozenge$
\end{thm}

\bigskip

\begin{pro} \cite{Diestel}, \cite{Pettis}\label{prop 5.1} If $x(.)$ is Pettis integrable and $h(.)$ is a measurable and essentially
bounded real-valued function, then $h(.)x(.)$ is Pettis integrable.$\hfill\lozenge$
\end{pro}

\bigskip

\begin{lem}\label{Cichon}\cite{Cichon}
Let $X$ be a Banach space, and $f : [a, b] \longrightarrow X$ be a Pettis integrable function.
Then
$$\displaystyle \int_a^b f(s)\,ds\in (b-a)\overline{co}\{f([a,b])\}.$$$\hfill\lozenge$
\end{lem}

\bigskip

\begin{thm}\label{app}
Suppose that  the conditions $(\bf C_0)$-$(\bf C_2)$ are satisfied. Then, the problem (\ref{NII})
 has a
solution $x$ in $\mathcal{B}_E(R)$ whenever 
 \begin{eqnarray}\label{max1}\left[\|q\|_\infty+\|k(\cdot,\cdot)\|_\infty \|h\|_{L_1}\,R\right]\|\vartheta_1\|_\infty R +2\|\vartheta_2\|_\infty R +\|T_2(0,x(0))\|\leq R\end{eqnarray}
and
\begin{eqnarray}\label{max2}[\|q\|_\infty+\|k(\cdot,\cdot)\|_\infty\,
\|h\|_{L^1}\,R]\,
\|
\vartheta_1\|_\infty\, r+2\,\left\|\vartheta_1\right\|_\infty\, R\, \,\left\|k(\cdot,\cdot)\right\|_\infty  \Phi(r) +\Psi(r)<r, \end{eqnarray}
 for $0<r\leq R.$
$\hfill\lozenge$\end{thm}
\noindent\textbf{Proof.}
Let $E=C(J,X)$ and
 $$S=\big\{x\in E:\, \|x\|_\infty\leq R\hbox{ and } \|x(t)-x(\tau)\|\leq b(t,\tau), \hbox{ for } t , \tau \in J\big\},$$
 where \begin{eqnarray*}b(t,\tau)&=&\frac{1}{(1-\delta_1\|\vartheta_1\|_\infty)}\Bigg
 \{R\delta_1
|\vartheta_1(t)-\vartheta_1(\tau)|+R\,|\vartheta_2(t)-\vartheta_2(\tau)|\\
&&+\delta_2 \big[\left\|q(t)-q(\tau)\right\|+R\sup_{s\in J}|k(t,s)-k(\tau,s))|\,\|h\|_{L^1}+ \displaystyle R\int_{t}^{\tau}|k(\tau,s)|\,h(s)\, ds\big]\Bigg\},\\
 \delta_1&=& \|q\|_\infty+\|k\|_\infty\|h\|_{L^1}R,\\
 \delta_2&=&\|\vartheta_1\|_\infty R+\|T_1(0,\theta)\|.\end{eqnarray*}
Obviously, $S$ is a closed, convex, bounded, and equi-continuous subset of $E.$
Let us consider three multi-valued mappings $A,B$ and $C$ defined on $S$ by:
\begin{eqnarray*}
  \begin{array}{rcl}&&\noindent (Ax)(t)=\displaystyle T_1(t,x(t)),\,\,
  \noindent (Cx)(t)=\displaystyle T_2(t,x(t)),\\
 &&\noindent Bx= \displaystyle \left\{q(t)+\int_0^{t}k(t,s)\,\xi(s)\, ds;\,\xi: J\longrightarrow X\hbox{ is P.I., } \xi(s)\in H(s,x(s)), a.e.\, s\in J\right\}.
\end{array}\end{eqnarray*}
The integral inclusion (\ref{in}) can be written in the form:
 \begin{eqnarray*}\label{eq r0}
\displaystyle x(t)\in Ax(t)\cdot Bx(t)+Cx(t).\end{eqnarray*}
In this section we will apply Theorem \ref{TT5}. So, the proof is composed by the following steps.
\par\vskip0.1cm
\noindent \textbf{Step 1:} $A$ and $C$ maps $S$ into $E$ and $B$ maps $S$ into $\mathcal{P}(E).$
\\
  The claim regarding $A$ and $C$ is clear in view of Condition $(\bf C_1).$ We corroborate now the claim for the multi-valued mapping $B.$
  Firstly, it follows from Conditions $(\bf C_2)$-$(\bf i),$ $(\bf C_2)$-$(\bf ii)$ and
  Proposition \ref{prop 5.1} in conjunction with Condition $(\bf C_0)$ that $B$ is well defined.
 Let $x\in S$ be arbitrary and let $y\in Bx.$ Then, there exists
a Pettis integrable mapping $\xi: J\longrightarrow X$ with $\xi(s)\in H(s,x(s))$ such that
$$y(t)=\displaystyle q(t)+\int_0^{t}k(t,s)\,\xi(s)\, ds, \, t\in J.$$
Let $\{t_n~,~ n\in \mathbb{N}\}$ be any sequence of elements in $J$ converging to $t.$ Then,  \begin{eqnarray*}\left\| y(t_n)-y(t)\right\|\leq\left\|q(t_n)-q(t)\right\|+\displaystyle \big\|\int_0^{t_n}(k(t_n,s)-k(t,s))\,\xi(s)\, ds\big\|+\big\|\int_{t_n}^{t}k(t,s)\,\xi(s)\, ds\big\|.
\end{eqnarray*}
From the Hahn-Banach theorem there exists $x^*\in X^*$ with
$\|x^*\|=1$ such that
$$\big\|\int_0^{t_n}(k(t_n,s)-k(t,s))\,\xi(s)\, ds\big\|=x^*\big(\int_0^{t_n}(k(t_n,s)-k(t,s))\,\xi(s)\, ds\big)$$
and
$$\big\|\displaystyle\int_{t_n}^{t} k(t,s)\xi(s)\, ds\big\|=x^*\big(\displaystyle\int_{t_n}^{t} k(t,s)\,\xi(s)\, ds\big)
.$$
Thus, in view of condition $(\bf C_2)$-$(\bf ii)$ we get \begin{eqnarray*}\| y(t_n)-y(t)\|&\leq&\|q(t_n)-q(t)\|+\int_0^{t_n}|k(t_n,s)-k(t,s))|\,\|\xi(s)\|\, ds+\int_{t_n}^{t}|k(t,s)|\,\|\xi(s)\|\, ds\\
&\leq&\|q(t_n)-q(t)\|+R\big(\int_0^{t_n}|k(t_n,s)-k(t,s))|\,h(s)\, ds+\int_{t_n}^{t}|k(t,s)|\,h(s)\, ds\big)\\
&:= &\|q(t_n)-q(t)\|+\sum_{i=1}^2I_i(t_n,t),
\end{eqnarray*}
where $$\displaystyle I_1(t_n,t)=R\int_0^{t_n}|k(t_n,s)-k(t,s))|\,h(s)\, ds\,\hbox{ and } \,\displaystyle I_2(t_n,t)=R\int_{t_n}^{t}|k(t,s)|\,h(s)\, ds.$$
Using condition $(\bf C_0)$ and  the fact that $h\in L^1(J,X)$ together with the dominated convergence theorem, we obtain $I_1(t_n,t)\rightarrow 0$ as $t_n\rightarrow t.$ Similarly, it is easy to see that $I_2(t_n,t)\rightarrow 0$ as $t_n\rightarrow t,$ since $h\in L^1(J,X)$ and $k\in L^\infty(J\times J).$
Hence, the continuity of $q$ achieves that $y(t_n) \rightarrow y(t),$ which implies that $B(x)\in \mathcal{P}(E).$
\par\vskip0.1cm
\noindent\textbf{Step 2:} $A$ is countably $\mathcal{D}$-set-lipschitzian on $S.$
  \\
 \noindent The use of condition $(\bf C_1)$-$(\bf ii)$ allows us to deduce that $A$ is lipschitzian with lipschitz constant $\|\vartheta_1\|_\infty$.
Since $T_1(t,\cdot)$ is sequential weak continuous, it is well-known that $A$ is $\mathcal{D}$-set-lipschitzian in light of Lemma $2.4$ in \cite{JKB2017} with  $\mathcal{D}$-function $\Phi_A(r)=\|\vartheta_1\|_\infty\, r,$ in particular it is countably $\|\vartheta_1\|_\infty$-set-lipschitzian.
\par\vskip0.1cm
\noindent\textbf{Step 3:} $B$ has a  weakly sequentially closed graph and is countably $\mathcal{D}$-set-lipschitzian.
\\
 \noindent To see this, let $(x_n)_n$ be a sequence of elements in $S$ weakly converging to some element $x\in S$ and let $y_n\in B(x_n)$ such that $y_n$ is weakly converging to some $y\in E.$ Thus, we can see that
there exists a  sequence of Pettis integrable mappings $(\xi_n)_n$ with $\xi_n(s)\in H(s,x_n(s))$   such that $$\displaystyle y_n(t)=q(t)+\int_0^{t}k(t,s)\,\xi_n(s)\, ds.$$
Fix $t\in I.$ Using the Dobrakov's Theorem, we result that $x_n(s)\rightharpoonup x(s)$ for each $s\in [0,t].$ Then the set $\{x_n(s), n\in \mathbb{N}\}$ is relatively weakly compact for each $s\in [0,t].$
 Using Condition $(\bf C_2)$-$(\bf iii)$ we obtain
\begin{eqnarray*}
\displaystyle \beta\left( H\left([0,t]\times \left\{x_n(s), n\in \mathbb{N}\right\} \right)\right)
&\leq &\displaystyle \beta \left(H\left(J\times \left\{x_n(s), n\in \mathbb{N}\right\} \right)\right)\\
&\leq &\displaystyle \Phi\left( \beta\left(\left\{x_n(s), n\in \mathbb{N}\right\} \right)\right).\end{eqnarray*}
Then, we get
\begin{eqnarray}\label{4}
\displaystyle  \beta\left( H\left([0,t]\times \left\{x_n(s), n\in \mathbb{N}\right\} \right)\right)=0.
\end{eqnarray}
This inequality in conjunction with inclusion
 $$\left\{\xi_n(s), n\in \mathbb{N}\right\}\subset H\left([0,t]\times \left\{x_n(s), n\in \mathbb{N}\right\} \right)\hbox{ for a.e. }s\in [0,t],$$
  allows us to infer that
  $\left\{\xi_n(s), n\in \mathbb{N}\right\}$ is relatively compact for a.e. $s\in [0,t].$ By the Eberlein-\v{S}mulian's Theorem,
 there exists a subsequence $(\xi_{n_k}(s))_k$ of $(\xi_n(s))_n$ which converges to some element $\xi(s).$
 Taking into account the fact that  $H(t,\cdot)$ has a weakly sequentially closed graph,  we get $\xi(s)\in H(s,x(s))$ a.e. $s\in J.$ This means that $$x^*\left(\xi_{n_k}(s)\right)\rightarrow x^*(\xi(s))\,\hbox{ for all } x^*\in
 X^*.$$
From the Condition  $(\bf C_2)$-$(\bf ii)$, it follows that $$\|\xi_n(s)\|\leq h(s)\, R.$$ The dominated convergence theorem for the Pettis integral reach to result that $\xi$ is Pettis integrable and
\begin{eqnarray}\label{y1}\displaystyle y(t)=q(t)+\int_0^{t}k(t,s)\,\xi(s)\, ds, \, t\in J.\end{eqnarray}
Next we prove that $B$ is countably $\mathcal{D}$-set-lipschitzian on $S$. Let $M$ be a countably bounded subset of $S.$ By the subadditivity property of $\beta$ and Lemma \ref{Cichon} together with Condition $(\bf C_2)$-$(\bf iii)$ we have for all $t\in J,$
 \begin{eqnarray*}\displaystyle
\beta(B(M)(t))&\leq&\beta\big(\bigcup\big\{\displaystyle\int_0^{t}k(t,s)\,H(s,x(s))\,ds,\, x\in M\big\}\big)\\
&\leq&\beta\big(\displaystyle\int_0^{t}k(t,s)\, H(s,M(s))\,ds\big)\\
&\leq&\|k(\cdot,\cdot)\|_\infty\beta\big(\displaystyle\overline{ co}(H(J,M(J))\big)\\
&\leq& \|k(\cdot,\cdot)\|_\infty\Phi(\beta(M(J))).
\end{eqnarray*}
A similar reasoning as in the
first steep, we can verify that $B(S)$ is an equi-continuous subset of $E$.
 Then invoking Theorem \ref{Mitchell} and the facts that $B(S)$ is bounded and $S$ is equi-continuous, we get $$\displaystyle
\beta(B(M))\leq \|k(\cdot,\cdot)\|_\infty\Phi(\beta(M)),$$ which implies that $B$ is a countably $\mathcal{D}$-set-lipschitzian with $\mathcal{D}$-function $\Phi_B(r)=\|k(\cdot,\cdot)\|_\infty\Phi(r).$
\par\vskip0.1cm
 \noindent\textbf{Step 4:} $C$ is weakly sequentially continuous and is countably $\mathcal{D}$-set-contraction.
\\
 \noindent From condition $(\bf C_1)$ it follows that $C$ is weakly sequentially   continuous.
Now,
let $M$ be a countably bounded subset of $S.$  By using condition $(\bf C_1)$-$(\bf iii)$ and the fact that $\Psi$ is nondecreasing, one sees that
 \begin{eqnarray*}\displaystyle
\beta(C(M)(t))
&\leq&\beta\big(\displaystyle T_2(t\times M(t))\big)\\
&\leq& \beta(\displaystyle(T_2(J\times M(t)))\\
&\leq& \Psi(\beta(M(t))).
\end{eqnarray*}
From Theorem \ref{Cichon} together with the boundedness of $S$ and the fact that $\Psi$ is nondecreasing, it follows that
\begin{eqnarray*}\displaystyle
\beta(C(M)(t))
&\leq& \Psi(\beta(M)).
\end{eqnarray*}
According to condition $(\bf C_1)$-$(\bf ii)$, it follows that $C(S)$ is an equi-continuous subset of $E.$
Taking into account the boundedness of $C(S),$ and using Theorem \ref{Mitchell} we obtain $$\displaystyle
\beta(C(M))\leq \Psi(\beta(M)).$$
Hence, $C$ define a countably $\mathcal{D}$-set-lipschitzian mapping on $S,$ with $\mathcal{D}$-function $\Phi_C=\Psi.$
\par\vskip0.1cm
\noindent\textbf{Step 5:} $(A\cdot B+C)(x)$ is a closed convex subset of $S$ for each $x\in S.$
\\
Let $x\in S$ and $y\in (A\cdot B+C)(x)$ such that
 \begin{eqnarray}\label{5}y(t)= T_1(t,x(t)) \big(q(t)+\int_0^{t}k(t,s)\zeta(s)\, ds\big)+T_2(t,x(t)),\end{eqnarray}
 where $\zeta : J \longrightarrow X$ is a Pettis integrable mapping with $\zeta(s)\in H(s,x(s)).$
From the Hann-Banach's theorem, it follows that \begin{eqnarray*}\|y(t)\|&\leq&  \left\|T_1(t,x(t))\right\|\, \big\|q(t)+\int_0^{t}k(t,s)\zeta(s)\, ds\big\|+\|T_2(t,x(t))\|\\
&\leq&  \|T_1(t,x(t))\| \big(\|q(t)\|+\int_0^{t}|k(t,s)|\,\|\zeta(s)\|\, ds\big)+\|T_2(t,x(t))\|\\
&\leq&  \|\vartheta_1\|_\infty R \left(\|q\|_\infty+\|k(t,\cdot)\|_\infty \|h\|_{L_1}\,R\right)+2\|\vartheta_2\|_\infty R +\|T_2(0,x(0))\|.
\end{eqnarray*}
Then, from our assumptions $\|y\|_\infty\leq R.$
Now, let $x\in S$ and $t, \tau \in J.$ Then, from equality (\ref{5}) it follows that
\begin{eqnarray*}\|y(t)-y(\tau)\|
&\leq&\big\|A(x(t)) \big(q(t)+\int_0^{t}k(t,s)\zeta(s)\, ds\big)-A(x(\tau))\big(q(\tau)+\int_0^{\tau}k(\tau,s)\zeta(s)\, ds\big)\big\|\\
&&+\big\|T_2(t,x(t))- T_2(\tau,x(\tau)) \big\|\\
&\leq&\big\| A(x(t))  \big(q(t)+\int_0^{t}k(t,s)\zeta(s)\, ds\big)-A(x(\tau)) \big(q(\tau)+\int_0^{\tau}k(\tau,s)\zeta(s)\, ds\big)\big\|\\
&&+\big\|T_2(t,x(t))- T_2(\tau,x(\tau)) \big\|\\
&\leq&\left\| A(x(t))-A(x(\tau))\right\| \big\| \big(q(t)+\int_0^{t}k(t,s)\zeta(s)\, ds\big)\big\|\\
&&+\left\|A(x(\tau)) \right\|\,\big\|\big(q(t)+\int_0^{t}k(t,s)\zeta(s)\, ds\big)-\big(q(\tau)+\int_0^{\tau}k(\tau,s)\zeta(s)\, ds\big)\big\|\\
&&
+\left\| C(x(t))-C(x(\tau))\right\|\\
&\leq &\sum_{i=1}^3 J_i(t,\tau).
\end{eqnarray*}
We have,
\begin{eqnarray*}
J_1(t,\tau)&=& \left\| A(x(t))-A(x(\tau))\right\| \big\| \big(q(t)+\int_0^{t}k(t,s)\zeta(s)\, ds\big)\big\|\\
&\leq&\delta_1\,\big(\vartheta_1(t)\,\|x(t)-x(\tau)\|
+|\vartheta_1(t)-\vartheta_1(\tau)| \,\|x(\tau)\|\big)\\
&\leq&\delta_1\,\big[\|\vartheta_1\|_\infty\,b(t,\tau)
+|\vartheta_1(t)-\vartheta_1(\tau)| \,R\big].
\end{eqnarray*}
Proceeding as in the steep $1,$ we obtain
\begin{eqnarray*}
J_2(t,\tau)&=&\left\|A(x(\tau)) \right\|\big\|\big(q(t)+\int_0^{t}k(t,s)\zeta(s)\, ds\big)-\big(q(\tau)+\int_0^{\tau}k(\tau,s)\zeta(s)\, ds\big)\big\|\\
&\leq&
\delta_2\,\big[\left\|q(t)-q(\tau)\right\|+R\sup_{s\in J}|k(t,s)-k(\tau,s))|\,\|h\|_{L^1}+ \displaystyle R\int_{t}^{\tau}|k(\tau,s)|\,h(s)\, ds\big].\end{eqnarray*}
Moreover, we have
\begin{eqnarray*}
J_3(t,\tau)&=&\left\|T_2(t,x(t))- T_2(\tau,x(\tau)) \right\|\\
&\leq& R|\vartheta_2(t)-\vartheta_2(\tau)|.\end{eqnarray*}
Accordingly, we have
\begin{eqnarray*}\|y(t)-y(\tau)\|
&\leq&\delta_1\big[\|\vartheta_1\|_\infty\,b(t,\tau)
+|\vartheta_1(t)-\vartheta_1(\tau)\ \,R\big]+R|\vartheta_2(t)-\vartheta_2(\tau)|\\
&&+\delta_2 \big[\left\|q(t)-q(\tau)\right\|+R\sup_{s\in J}|k(t,s)-k(\tau,s))|\,\|h\|_{L^1}+ \displaystyle R\int_{t}^{\tau}|k(\tau,s)|\,h(s)\, ds\big].\end{eqnarray*}
By using the equality,
\begin{eqnarray*}(1-\delta_1\|\vartheta_1\|_\infty)\,b(t,\tau)&=&
\delta_1\,|\vartheta_1(t)-\vartheta_1(\tau)| \,R+R\,|\vartheta_2(t)-\vartheta_2(\tau)|
+\delta_2 \big[\left\|q(t)-q(\tau)\right\|\\
&&+R\sup_{s\in J}|k(t,s)-k(\tau,s))|\,\|h\|_{L^1}+ \displaystyle R\int_{t}^{\tau}|k(\tau,s)|\,h(s)\, ds\big],\end{eqnarray*}
we obtain
\begin{eqnarray*}\|y(t)-y(\tau)\|
\leq b(t,\tau).\end{eqnarray*}
Consequently, we have $$(A\cdot B+C)(x)\in \mathcal{P}(S).$$
 Next we will prove that $A\cdot B+C$ has convex values.
 Let $x\in S$ and let $y_1, y_2\in Ax\cdot Bx+Cx.$
Then,
 $$y_1(t)=T_1(t,x(t)) \big(q(t)+\int_0^{t}k(t,s)\zeta_1(s)\, ds\big)+ T_2(t,x(t)),$$
and
$$y_2(t)=T_1(t,x(t)) \big(q(t)+\int_0^{t}k(t,s)\zeta_2(s)\, ds\big)+T_2(t,x(t)),$$
 where $\zeta_1, \zeta_2: J \longrightarrow X$ are two Pettis integrable mappings with $\zeta_1(s),\zeta_2(s)\in H(s,x(s)).$
So, for any $\alpha \in (0,1),$ we have
\begin{eqnarray*} \alpha y_1(t)+(1-\alpha)y_2(t)
&=& \displaystyle  \alpha T_1(t,x(t))\big(q(t)+\int_0^{t}k(t,s)\zeta_1(s)\, ds\big)\\
&&+(1-\alpha)T_1(t,x(t))\big(q(t)+\int_0^{t}k(t,s)\zeta_2(s)\, ds\big)+T_2(t,x(t)) \\
&=& \displaystyle T_1(t,x(t))\,q(t)+T_1(t,x(t))\,\int_0^{t}k(t,s)(\alpha\zeta_1(s)+(1-\alpha)\zeta_2)\, ds+T_2(t,x(t)).
\end{eqnarray*}
Since $H$ has convex values, we get $$\alpha\zeta_1(s)+(1-\alpha)\zeta_2(s)\in H(s,x(s)) \hbox{ a.e. }s\in J.$$
From the above discussion, it is clear that  $ \alpha\zeta_1(\cdot)+(1-\alpha)\zeta_2(\cdot)$ is Pettis integrable and we have $$\alpha y_1+(1-\alpha)y_2\in Ax \cdot Bx+Cx.$$
It remains to prove that $Ax \cdot Bx+Cx$ is closed. Let $y_n$ be a sequence of elements in $Ax \cdot Bx+Cx$ which converging to some $y\in E.$
Then, $$y_n(t)=T_1(t,x(t)) \big(q(t)+\int_0^{t}k(t,s)\xi_n(s)\, ds\big)+ T_2(t,x(t)),$$
where
 $\xi_n: J \longrightarrow X$ is a sequence of Pettis integrable mapping with  $\xi_n(s)\in H(s,x(s)).$ The use of Condition $(\bf C_2)$-$(\bf iii),$ yields that $H$ has relatively weakly compact values. Thus, by using the inclusion $$(\xi_n(s))_n\subset H(s,x(s)),$$
it follows that $(\xi_n(s))_n$ has a convergent subsequence $(\xi_{n_k}(s))_k$ to some $\xi\in H(s,x(s)),$ in view of the Eberlein-\v{S}mulian's Theorem. From our assumptions and by using the dominated convergence theorem of Pettis we  can deduce that $\xi$ is Pettis integrable and $$\displaystyle\int_0^{t}k(t,s)\xi_n(s)\, ds\rightarrow \int_0^{t}k(t,s)\xi(s)\, ds.$$
This implies that $$y_n(t)\rightarrow T_1(t,x(t)) \big(q(t)+\int_0^{t}k(t,s)\xi(s)\, ds\big)+ T_2(t,x(t)) ,$$ and consequently $y\in Ax \cdot Bx+Cx.$
\par\vskip0.1cm
\noindent\textbf{Step 6:} $\|B(S)\|\,\Phi_A(r)+\|A(S)\|\,\Phi_B(r)+\Phi_A(r)
\,\Phi_B(r)+\Phi_C(r)<r$ for $0<r<\|S\|.$\\
First let's note that $A(S)$ and $B(S)$ are bounded subsets with bounds  $\|\vartheta_1\|_\infty\, \|S\|$ and $\|q\|_\infty+\|k(\cdot,\cdot)\|_\infty\,\|h\|_{L^1}\|S\|.$ From the above steeps, we have $\Phi_A(r)=\|\vartheta_1\|_\infty\, r,$ $\Phi_B(r)=\|k(\cdot,\cdot)\|_\infty\, \Phi(r)$ and $\Phi_C(r)=\Psi(r),$ respectively. From the estimate (\ref{max2}) and since $S\subset \mathcal{B}_E(R)$ we obtain $$\left[\|q\|_\infty+\|k(\cdot,\cdot)\|_\infty\,\|h\|_{L^1}\|S\|\right]\|\,\vartheta_1\|_\infty\, r+\|\vartheta_1\|_\infty\, \|S\|\,\|k(\cdot,\cdot)\|_\infty\, \Phi(r)+\|\vartheta_1\|_\infty\,\|k(\cdot,\cdot)\|_\infty\, r\,\Phi(r) +\Psi(r)<r,$$ for all $0<r\leq \|S\|.$ $\hfill\Box$

\bigskip

 As an application, we discuss existence results for the following nonlinear functional differential inclusion (FDI):
\begin{eqnarray}\label{in2}
 \displaystyle
 \displaystyle \left(\displaystyle\frac{ x(t)-T_2(t,x(t))}{T_1(t,x(t))}\right)'\in k(t)\,H(t,x(t)), \, t\in J,
\end{eqnarray}
satisfying the initial condition
\begin{eqnarray}\label{in3}x(0)=x_0\in X,\end{eqnarray}
where $J=[0,1],$ $k :  J \longrightarrow \mathbb{R}$, $T_i : J \times X \longrightarrow X$ and $ H : J \times X \longrightarrow \mathcal{P}_{cp,cv}(X).$
 By a solution of the problem (\ref{in2})-(\ref{in3}) we mean a function $x\in C(J,X)$
 that satisfies the inclusions (\ref{in2})-(\ref{in3}) on $J.$
 Our existence result for  the (FDI) (\ref{in2})-(\ref{in3}) is
\begin{thm}\label{app}
Suppose that the conditions $(\bf C_1)$-$(\bf C_3)$ hold.
Further if there exists a real number $R > 0$ such that the inequalities (\ref{max1}) and (\ref{max2}) hold
 with $\|q\|_\infty=\displaystyle\frac{\|x_0-T_2(0,x_0)\|}{\|T_1(0,x_0)\|},$ then the (FDI) (\ref{in2})-(\ref{in3})
 has a
solution $x$ in $\mathcal{B}_R.$~~$\hfill\lozenge$\end{thm}
\noindent\textbf{Proof.}
Notice first that the (FDI) (\ref{in2})-(\ref{in3})
can be written as a fixed point problem
\begin{eqnarray*}\label{in}
 \displaystyle
\left\{
  \begin{array}{lll}
 x(t)\in T_1(t,x(t))\cdot \left[\displaystyle\frac{x_0-T_2(0,x_0)}{T_1(0,x_0)}+ \displaystyle\int_0^{t}k(s)\,H(s,x(s))\, ds\right]+T_2(t,x(t)),\, t\in J, \\
  \displaystyle  x(0)= x_0.
  \end{array}
\right.
\end{eqnarray*}
Now the desired conclusion of the theorem follows by a direct application of
Theorem \ref{app} with $\displaystyle q(t) = \frac{x_0-T_2(0,x_0)}{T_1(0,x_0)}$ and $k(t,s)=k(s)$
for all $t,s \in J.$ $\hfill\Box$
\bigskip
\begin{rem} Finally while concluding this work, we remark if $T_1(t,x) = 1_X$ on
$J \times X$ and $q(t)=\theta$ on $J,$ then Theorem \ref{app} reduces to the existence results for the Perturbed Volterra Integral inclusion
\begin{eqnarray*}\label{in}
 \displaystyle
\displaystyle x(t)\in \int_0^{t}k(t,s)\,H(s,x(s))\, ds+T(t,x(t)), \, t\in J
 .
\end{eqnarray*}$\hfill\lozenge$
\end{rem}

\end{document}